\documentclass[reqno]{amsart}
\usepackage{amssymb,latexsym,amscd}
\theoremstyle{plain}
\newtheorem{theorem}{Theorem}
\newtheorem{lemma}{Lemma}
\newtheorem{proposition}{Proposition}
\newtheorem{corollary}{Corollary}
\theoremstyle{definition}
\newtheorem{definition}{Definition}
\newtheorem{remark}{Remark}

\newcommand{\mg}{\mathfrak{g}}
\newcommand{\mgs}{\mathfrak{g^{*}}}
\newcommand{\f}{\mathfrak}
\newcommand{\bhat}[1]{\hat{\textbf{#1}}}
\newcommand{\bb}[1]{\textbf{#1}}
\begin{document}
\begin{abstract}
In this article\footnote{} we study conditions for the integrability
of the distribution defined on a regular Poisson manifold as the
orthogonal complement (with respect to a pseudo-Riemannian metric)
to the tangent spaces of the leaves of a symplectic foliation.
Examples of integrability and non-integrability of this
distribution are provided.
\end{abstract}
\title[On the integrability]%
{On the integrability of orthogonal distributions in Poisson manifolds}
\author{Daniel Fish}
\address{Department of Mathematics and Statistics, Portland State University,
Portland, OR, U.S.} \email{djf@pdx.edu}
\author{Serge Preston}
\address{Department of Mathematics and Statistics, Portland State University,
Portland, OR, U.S.} \email{serge@mth.pdx.edu}
 \maketitle

 \section{Introduction}


Let $(M^{n},P)$ be a regular Poisson manifold. Denote by ${\mathcal S}=\{S_{m}\vert m\in M\}$
the symplectic foliation of $M$ by symplectic leaves (of constant dimension $2\leq k <n$ in the regular case).
 Denote by $T(\mathcal S)$ the sub-bundle of $T(M)$ of tangent spaces to the symplectic leaves
 (the association $x\rightarrow T_{x}(\mathcal S)$ is an integrable distribution on $M$ which we
  will also denote by $T(\mathcal S)$). Let $M$ be endowed with a pseudo-Riemannian metric $g$
  such that the restriction of $g$ to each symplectic leaf is non-degenerate By continuity, the signature of the restriction of $g$ to $T_{m}({\mathcal S})$ is the same for all $m\in M$.

\par Let $\mathcal{N}_{m}=S_{m}^{\bot}$ be the subspace of $T_{m}(M)$ that is $g$-orthogonal to $S_{m}$. The
 association $m\rightarrow \mathcal{N}_{m}$ defines a distribution $\mathcal{N}$ which is transversal
 and complemental to the distribution $T(\mathcal S)$. The restriction of the metric $g$ to $\mathcal N$
  is non-degenerate and has constant signature. In general, the distribution $\mathcal N$ is not integrable.

\par If the metric $g$ is Riemannian, and if the Poisson tensor is parallel with respect to the Levi-Civita
 connection $\nabla =\nabla^{g}$ defined by $g$, ie: $\nabla P=0$, then it is a classical result of A. Lichnerowicz
 (see Vaisman 1994, Remark 3.11) that the distribution $\mathcal{N}$ is integrable,
  and the restriction of the metric $g$ to the symplectic leaves defines, together with the symplectic
   structure $\omega_{S}=P\vert_{S}^{-1}$, a {\bf K\"{a}hler structure} on symplectic leaves.
\par

Integrability of the distribution $\mathcal N$ depends strongly on
the foliation $\mathcal S$ and its ``transversal topology" (see
Molino 1988, Reinhart 1983, Ch.4 for the Riemannian case). Thus,
in general it is more a question of the theory of bundles with
Ehresmann connections rather than that of Poisson geometry. Yet in
some instances, it is useful to have integrability conditions in
terms of the Poisson structure $P$, and to relate integrability of
the distribution $\mathcal N$ with other structures of the Poisson
manifold - Casimir functions, Poisson vector fields, etc.
\par

Our interest in this question was influenced by our study of the
representation of a dynamical system in metriplectic form, i.e. as
a sum of a Hamiltonian vector field (with respect to a Poisson
structure, (see Bloch et al 1996, Morrison 1986, Grmela 1990, Fish
2005) and a gradient one (with respect to a metric $g$).
Integrability of the distribution $\mathcal N$  guarantees that in
the geometrical (local) splitting of the space $M$ as a product of
a symplectic leaf and a transversal submanifold with Casimir functions $c^{i}$ as local coordinates (see Weinstein 1983), the
transversal submanifold can be chosen to be invariant under the
gradient flow (with respect to the metric $g$) of the functions $c^{i}$.\par

As a result one can separate observables of the system into
Casimirs undergoing pure gradient (dissipative) evolution from
those (along symplectic leaves) which undergo the mix of
Hamiltonian and gradient evolutions.  Such a separation leads to
an essential simplification of the description of the transversal
dynamics in metriplectic systems.\par

The structure of this article is as follows. In Section 2 we introduce necessary notions and notations. In Section 3 we obtain necessary and sufficient conditions on the
metric $g$ and the tensor $P$ for the distribution $\mathcal{N}$
to be integrable. We derive these conditions in terms of covariant derivatives of the Poisson
Tensor, in terms of covariant derivatives of Casimir covectors, and
as conditions on the nullity of the Nijenhuis Torsion of the
(1,1)-tensor $A^{\mu }_{\nu }=P^{\mu \kappa }g_{\kappa \nu}$. As a
corollary we prove that the distribution $\mathcal N$ is
integrable if parallel translation (via the Levi-Civita connection of the metric $g$) in the direction of $\mathcal N$
preserves the symplectic distribution $T(\mathcal S)$.

\par In Section 4 we present integrability conditions in Darboux-Weinstein coordinates: the distribution
$\mathcal N$ is integrable if and only if the following symmetry conditions are fulfilled for
$\Gamma$
\[
\Gamma_{JIs}=\Gamma_{JIs},
\]
where $\Gamma_{\alpha \beta \gamma}=g_{\alpha \sigma}\Gamma^{\sigma}_{ \beta \gamma}$, and where
capital Latin letters $I,J$ indicate the transversal coordinates ,while small Latin letters indicate
coordinates along symplectic leaves.\par

In Section 5 we describe some examples of non-integrability, describe a model example of a 4d Poisson
 manifold with Poisson structure of rank $2$, where the distribution $\mathcal{N}$ is not integrable, and discuss nonintegrability in the case of a topologically nontrivial symplectic fibration.
\par
In Section 6 we prove integrability of $\mathcal N$ for linear
Poisson structures on dual spaces $\mgs$ of real semi-simple Lie
algebras $\mg$, with the metric $g$ induced by the Killing form,
as well as for the dual $\mathfrak{se}(3)^{*}$ to the Lie algebra $\mathfrak{se}(3)$ of
Euclidian motions with the simplest non-degenerate
$ad^{*}$-invariant metric(s) (see Zefran et al. 1996).

\section{Orthogonal distribution of Poisson manifold with a pseudo-Riemannian metric}
Let $(M^{n},P)$ be a regular Poisson manifold. We will be use
local coordinates $x^{\alpha }$ in the domains $U\subset M$ with
the corresponding local frame $\{\frac{\partial}{\partial
x^{\alpha}} \}$ and the dual coframe ${dx^{\alpha}}$. Let $g$ be a
pseudo-Riemannian metric on $M$ as above, and let $\Gamma$ denote
the Levi-Civita connection associated with $g$. The tensor
$P^{\tau \sigma }(x)$ defines a mapping
\[
0 \rightarrow C(M) \rightarrow T^{*}(M) \stackrel{P}{\rightarrow } T(\mathcal S)\rightarrow 0
\]
where $C(M)\subset T^{*}(M)$ is the kernel of $P$  and $T(\mathcal
S)$ is (as defined above) the tangent distribution of the
symplectic foliation $\{\mathcal{S}^{k} \}$. The space $C(M)$ is a
sub-bundle of the cotangent bundle $T^{*}M$ consisting of Casimir
covectors. Locally, $C(M)$ is generated by differentials of
functionally independent Casimir functions $c^{i}(x),\ i=1,\ldots
, n-k$ satisfying the condition $P^{\tau \sigma }dc^{i}_{\sigma }
=0$ (in this paper we assume the condition of summation by repeated
indices).
\par
We denote by $\mathcal{N}$ the {\bf distribution defined as the
$g$-orthogonal complement $T(\mathcal S)^{\perp}$ to $T(\mathcal
S)$ in $T(M)$}. Then we have, at every point $x$ a decomposition
into a direct sum of distributions (sub-bundles)
\[
T_{x}M = T_{x}(\mathcal S) \oplus \mathcal{N}_{x}.
\]

The assignment $x \rightarrow \mathcal{N}_{x}$ defines a {\bf
transverse connection} for the foliation $\mathcal S$, or, more
exactly, for the bundle $(M,\pi, M/\mathcal S)$ over the {\bf
space of leaves} $M/\mathcal S$, whenever one is defined (see
below). We are interested in finding necessary and sufficient
conditions on $P$ and $g$ under which the distribution
$\mathcal{N}$ is integrable.  By the Frobenius theorem,
integrability of $\mathcal{N}$ is equivalent to the involutivity
of the distribution $\mathcal{N}$ with respect to the Lie bracket
of $\mathcal{N}$-valued vector fields (sections of the sub-bundle
$\mathcal{N}\subset T(M)$).
\par Let $\omega ^{i} = \omega^{i}_{\mu }dx^{\mu }$ ($i\leq d=n$-$k$)
be a local basis for $C(M)$. For any $\alpha $ in $T^{*}(M)$, let
$\alpha^{\sharp}$ denote the image of $\alpha$ under the bundle
isomorphism $\sharp:T^{*}(M) \rightarrow TM$ of index lifting
induced by the metric $g$. The inverse isomorphism (index lowering)
will be conventionally denoted by $\flat :T(M) \rightarrow
T^{*}M$. We introduce the following vector fields: $\xi_{i}\,=(\omega^{i})^{\sharp } \in T(M)$.

\begin{lemma} The vectors $\xi_{i}$ form a (local) basis for $\mathcal{N}$. \end{lemma}
\begin{proof} Since $g$ is non-degenerate, the vectors $\xi_{i}$ are linearly independent and span
 a subspace of $TM$ of dimension $d$. For any vector $\eta \in TM$,
   \[<\xi_{i},\eta >_{g} \, = \, g_{\mu \nu}\xi_{i}^{\mu}\eta^{\nu} \, = \, g_{\mu \nu}g^{\mu
    \lambda}\omega_{\lambda}^{i}\eta^{\nu} \, = \, \omega_{\nu}^{i}\eta^{\nu} \, = \,  \omega^{i}(\eta ).\]
So the vector $\eta $ is $g$-orthogonal to all $\xi_{i}$ if and only if $\eta $ is
 annihilated by each $\omega ^{i}$. That is, $\, \eta \in Ann(C(M))= \{\lambda \in  T(M)\, |\,
 \omega^{j}(\lambda ) = 0 , \,\,\, \forall \, j \leq d \}$. Since $Ann(C(M)) = T(S)$,
  we see that the linear span of $\{\xi_{i}\}^{\perp} $ is $ T(\mathcal S)$.
\end{proof}

\begin{definition} The \textbf{curvature} (Frobenius Tensor) of the ``transversal connection''
 $\mathcal{N}$ is defined as the bilinear
mapping
                   \[
                   \mathfrak{R}_{\mathcal N}:T(M) \times T(M) \rightarrow T(\mathcal S)
                   \]
defined by
               \begin{equation}
               \mathfrak{R}_{\mathcal N}(\gamma , \eta) = v([h\gamma,h\eta ]),
               \end{equation}
where $h:T(M) \rightarrow \mathcal{N}$ is $g$-orthogonal projection onto $\mathcal{N}$, and $v:T(M)
\rightarrow T(\mathcal S)$ is $g$-orthogonal projection onto $T(\mathcal S)$.
\end{definition}
It is known (see DeLeon, Rodrigues 1989, Sec.1.15) that
$\mathcal{N}$ is integrable if and only if the curvature
$\mathfrak{R}_{\mathcal N}$ defined above is identically zero on
$TM \times TM$.
\begin{remark}
An equivalent way to characterize the integrability of
$\mathcal N$ is to use the structure tensor of J. Martinet or the
D. Bernard structure tensor of the annihilator ${\mathcal
N}^{*}\subset T^{*}(M)$ of the distribution $\mathcal N$ (see
Libermann 1976).
\end{remark}
\section{Integrability criteria.}
The condition (1) is equivalent to $v([\gamma ,\eta ] ) = 0$, for all $\gamma ,\eta \in \mathcal{N}$. If we write the vectors $\gamma , \eta$ in terms of the basis $\{\xi_{i}\}$, then we have
\begin{eqnarray*}
v([\gamma^{i}\xi_{i},\eta^{j}\xi_{j}]) & = & v\left(\gamma^{i}(\xi_{i}\cdot \eta^{j})\xi_{j}-\eta^{j}
(\xi_{j}\cdot \gamma^{i})\xi_{i} + \gamma^{i}\eta^{j}[\xi_{i},\xi_{j}]\right)\\
   & = & \gamma^{i}\eta^{j}v([\xi_{i},\xi_{j}]),\,\,\,\,\textrm{since}\,\,\,v(\xi_{k})=0\,\,
 \forall \, k.
\end{eqnarray*}
Thus $\mathfrak{R}=0$ if and only if $v([\xi_{i},\xi_{j}])=0$ \,\, for all $i,j \leq d$.

Consider the linear operator $A:T(M) \rightarrow T(M)$ defined by
the (1,1)-tensor field $A^{\tau}_{\mu} = P^{\tau \sigma }g_{\sigma
\mu}$. Since $g$ is non-degenerate we have $ImA = T(\mathcal S)$.
Since each basis vector $\xi_{i} \in \mathcal{N}$ is of the form
$\xi_{i}^{\mu} = g^{\mu \nu }\omega^{i}_{\nu}$ with $\omega^{i}
\in kerP$, we also have
 \[
 A^{\tau}_{\mu}\xi_{i}^{\mu}\, = \, P^{\tau \sigma}g_{\sigma \mu}g^{\mu \nu}\omega^{i}_{\nu}\,  = \,
  P^{\tau \nu}\omega^{i}_{\nu} = 0.
 \]
Therefore $\mathcal{N} \subset kerA$, and by comparing dimensions we see that $\mathcal{N} = kerA$. Notice that operator $A$ and the orthonormal projector $v$ have the same image and kernel. We conclude that
\[
\mathfrak{R} = 0 \, \Leftrightarrow \, A[\xi_{i},\xi_{j}]=0, \,\,\, \forall \, i,j \leq d.
\]
We now prove the main result of this section.

\begin{theorem} Let $\omega ^{i},\ 0\leq i\leq d$ be a local basis for $C(M)$
 and let $(\omega^{i})^{\sharp} = \xi_{i}$ be the corresponding local basis of $\mathcal N$.
  Let $\nabla$ be the Levi-Civita covariant derivative on $TM$ corresponding to the metric $g$.
   Then the following statements are equivalent:
\begin{enumerate}
\item The distribution $\mathcal{N}$ is integrable.
\item For all $ i,j \leq d$, and all $\tau \leq n$,
    \begin{equation}
  P^{\tau \sigma}(\nabla_{\xi_{i}}\omega^{j}_{\sigma}-\nabla_{\xi_{j}}\omega^{i}_{\sigma})= 0.
    \end{equation}
\item  For all $ i,j \leq d$, and all $\tau \leq n$,
    \begin{equation} g^{\lambda \alpha}(\nabla_{\lambda}P)^{\tau \sigma}
    (\omega^{i}\wedge \omega^{j})_{\sigma \alpha}= 0,
    \end{equation}
    where $\nabla_{\lambda}=\nabla_{\partial/\partial x^{\lambda}}$.
\item For all $ i,j \leq d$, and all $\tau \leq n$,
    \begin{equation} P^{\tau \sigma}g_{\sigma \lambda}(\nabla_{\xi_{i}}\xi_{j}^{\lambda} -
\nabla_{\xi_{j}}\xi_{i}^{\lambda}) = 0.
\end{equation}
\item The sub-bundle $C(M)$ is invariant under the following skew-symmetric bracket on 1-forms generated by the bracket of vector fields:
\begin{equation}
[\alpha ,\beta ]_{g}=[\alpha^{\sharp },\beta^{\sharp }]^{\flat }
\end{equation}
i.e. if \, $\alpha ,\beta \in \Gamma(C(M))$, then $[\alpha ,\beta ]_{g}\in \Gamma(C(M)).$
\end{enumerate}
\end{theorem}
\begin{proof}Since the Levi-Civita connection of $g$ is torsion-free, we know that
      \[ [\xi_{i},\xi_{j}] = \nabla_{\xi_{i}}\xi_{j} - \nabla_{\xi_{j}}\xi_{i}.\]
Therefore, in a local chart $(x^{\alpha})$,
\begin{eqnarray}
\nonumber A^{\tau}_{\lambda}[\xi_{i},\xi_{j}]^{\lambda} & =
& A^{\tau}_{\lambda}(\nabla_{\xi_{i}}\xi_{j}^{\lambda}-\nabla_{\xi_{j}}\xi_{i}^{\lambda})\\
\nonumber & = & P^{\tau \sigma}g_{\sigma \lambda}(\nabla_{\xi_{i}}\xi_{j} - \nabla_{\xi_{j}}\xi_{i})\\
\nonumber & = & P^{\tau \sigma}(\nabla_{\xi_{i}}\omega^{j}_{\sigma} -
\nabla_{\xi_{j}}\omega^{i}_{\sigma}).
\end{eqnarray}
In the last step we have used the fact that lifting and lowering of
indices by the metric $g$ commutes with the covariant derivative
$\nabla $ defined by the Levi-Civita connection of $g$.\par

 Recalling from the discussion before the Theorem that the
 integrability of the distribution $\mathcal{N}$ is equivalent to the
 nullity of $A[\xi_{i},\xi_{j}]$ for all $i,j\leqq d$, we see that statements (1), (2), and (4) are equivalent.
          \par To prove the equivalence of these statements to (3) we notice that
\begin{eqnarray}
\nonumber P^{\tau \sigma}(\nabla_{\xi_{i}}\omega^{j}_{\sigma} - \nabla_{\xi_{j}}\omega^{i}_{\sigma}) &=& P^{\tau \sigma}(\xi_{i}^{\lambda}\nabla_{\lambda}\omega^{j}_{\sigma} - \xi_{j}^{\lambda} \nabla_{\lambda}\omega^{i}_{\sigma})\\
\nonumber & = & P^{\tau \sigma}g^{\lambda \alpha}(\omega^{i}_{\alpha}\nabla_{\lambda}\omega^{j}_{\sigma}-\omega^{j}_{\alpha} \nabla_{\lambda}\omega^{i}_{\sigma})\\
\nonumber & = &  g^{\lambda \alpha}(\nabla_{\lambda}P)^{\tau \sigma}(\omega^{i}_{\sigma}\omega^{j}_{\alpha}- \omega^{i}_{\alpha}\omega^{j}_{\sigma})\\
\nonumber & = &  g^{\lambda \alpha}(\nabla_{\lambda}P)^{\tau
\sigma}(\omega^{i} \wedge \omega^{j})_{\sigma \alpha}.
\end{eqnarray}
Here, at the third step we have used the following equality
\[
P^{\tau \sigma}g^{\lambda
\alpha}\omega^{i}_{\alpha}\nabla_{\lambda}\omega^{j}_{\sigma}=-g^{\lambda
\alpha}(\nabla_{\lambda}P)^{\tau
\sigma}\omega^{i}_{\sigma}\omega^{j}_{\alpha},
\]
since $ P^{\tau \sigma}\omega^{j}_{\sigma}=0$ (similarly for the second term).
\par
To prove the equivalence the (5) with the other statements we act
as follows. Let $\alpha = \alpha_{i}\omega^{i}$ and $\beta =
\beta_{j}\omega^{j}$ be any two sections of sub-bundle $C(M)\subset
T^{*}(M)$. Then $\alpha^{\sharp}=\sum_{i}\alpha_{i}\xi_{i}$ and
$\beta^{\sharp} = \sum_{i}\beta_{j}\xi_{j}$. So we have
\begin{align}
\notag [\alpha,\beta ]_{g}   &=
\nabla_{\beta_{j}\xi_{j}}(\alpha_{i}\omega^{i})-
\nabla_{\alpha_{i}\xi_{i}}(\beta_{j}\omega^{j}) \\ \notag
&= \beta_{j}\left[\alpha_{i}\nabla_{\xi_{j}}\omega^{i} +
\frac{\partial \alpha_{i}}{\partial
x^{k}}\xi_{j}^{k}\omega^{i}\right]
 -\alpha_{i}\left[\beta_{j}\nabla_{\xi_{i}}\omega^{j}+
 \frac{\partial \beta_{j}}{\partial x^{k}}\xi_{i}^{k}\omega^{j}\right]\\ \notag
&=  \alpha_{i}\beta_{j}(\nabla_{\xi_{j}}\omega^{i} -
\nabla_{\xi_{i}}\omega^{j} ) + \beta^{\sharp}(\alpha_{i})\omega^{i} -
 \alpha^{\sharp}(\beta_{j})\omega^{j}\\ \notag
&=  \alpha_{i}\beta_{j}(\nabla_{\xi_{j}}\omega^{i} -
\nabla_{\xi_{i}}\omega^{j} ) +  (\sum_{j}\beta_{j}\xi_{j}^{k})
\frac{\partial \alpha_{i}}{\partial x^{k}}\omega^{i} -
(\sum_{i}\alpha_{i}\xi_{i}^{k})
\frac{\partial \beta_{j}}{\partial x^{k}}\omega^{j}\\
&=  \alpha_{i}\beta_{j}[\omega^{i},\omega^{j}]_{g} +
(\beta^{\sharp}(\alpha_{i}) -
\alpha^{\sharp}(\beta_{i}))\omega^{i}.
\end{align}
At the last step we have used the following (recall that
$\omega^{i}=(\xi_{i})^{\flat}$)
\[
\nabla_{\xi_{j}}\omega^{i} - \nabla_{\xi_{i}}\omega^{j}
=(\nabla_{\xi_{j}}\xi_{i} -
\nabla_{\xi_{i}}\xi_{j})^{\flat}=[\xi_{i},\xi_{j}]^{\flat}=[\omega^{i\
\sharp},\omega^{j\ \sharp}]^{\flat}=[\omega^{i},\omega^{j}]_{g}.
\]
The second term in (6) is always in the kernel
$C(M)$ of $P$, thus applying $P$ to both sides yields
\[P^{\tau \sigma}([\alpha,\beta]_{g})_{\sigma} = \, \alpha_{i}\beta_{j}P^{\tau \sigma}
([\omega^{i},\omega^{j}]_{g})_{ \sigma}=\,
 \alpha_{i}\beta_{j}P^{\tau \sigma}g_{\sigma \lambda}(\nabla_{\xi_{i}}\xi_{j}^{\lambda} -
\nabla_{\xi_{j}}\xi_{i}^{\lambda}).\] Therefore, condition (4)
above holds if and only if the space of sections of the bundle
$C(M)$ of Casimir covectors is invariant under the bracket
$[-,-]_{g}$.
\end{proof}
\vskip0.4cm

\begin{corollary} If $(\nabla_{(\omega^{i}) ^{\sharp}}P)^{\tau \sigma}\omega^{j}_{\sigma}=0$ for
 all $\sigma$, $i$ and $j$, i.e. if \, $\nabla_{(\omega^{i}) ^{\sharp}}P\vert_{C(M)}=0$ for all $i$,
 then the distribution $\mathcal N$ is integrable.
\end{corollary}
\begin{proof} In the proof of the equivalence of statements (1) and (2) with statement (3) in the Theorem, it was shown that
\begin{align*}
P^{\tau \sigma}(\nabla_{\xi_{i}}\omega^{j}_{\sigma} -
\nabla_{\xi_{j}}\omega^{i}_{\sigma}) &=  g^{\lambda
\alpha}(\nabla_{\lambda}P)^{\tau
\sigma}(\omega^{i}_{\alpha}\omega^{j}_{\sigma}-
\omega^{j}_{\alpha}\omega^{i}_{\sigma})\\
&=g^{\lambda
\alpha}\omega^{i}_{\alpha}(\nabla_{\lambda}P)^{\tau
\sigma}\omega^{j}_{\sigma}-g^{\lambda
\alpha}\omega^{j}_{\alpha}(\nabla_{\lambda}P)^{\tau
\sigma}\omega^{i}_{\sigma}\\
&=\xi^{\lambda}_{i}(\nabla_{\lambda}P)^{\tau
\sigma}\omega^{j}_{\sigma}-\xi_{j}^{\lambda}(\nabla_{\lambda}P)^{\tau
\sigma}\omega^{i}_{\sigma}\\
&=(\nabla_{\xi_{i}}P)^{\tau
\sigma}\omega^{j}_{\sigma}-(\nabla_{\xi_{j}}P)^{\tau
\sigma}\omega^{i}_{\sigma}.
\end{align*}
Since $\xi_{i}=(\omega^{i}) ^{\sharp}$ for each $i$, if $\, \nabla_{(\omega^{i}) ^{\sharp}}P\vert_{C(M)}=0$, then condition (3) of the
Theorem is fulfilled.
\end{proof}
\noindent The following criteria specify the part of the A.
Lichnerowicz condition that $P$ is $g$-parallel (see Vaisman 1994)
ensuring the integrability of the distribution $\mathcal N$:
\begin{corollary} If \, $\nabla_{\alpha^{\sharp}}:T(M) \rightarrow T(M)$ preserves the
tangent sub-bundle $T(\mathcal S)$ to the symplectic leaves for every $\alpha \in C(M)$, then $\mathcal{N}$ is
integrable.
\end{corollary}
\begin{proof} If the parallel translation $\nabla_{\alpha^{\sharp}}$ along the trajectories of the vector field
 $\xi = \alpha^{\sharp} $ preserves $T(\mathcal S)$, then it also
preserves its $g$-orthogonal complement $\mathcal{N}$, and hence
the dual to parallel translation in the cotangent bundle will
preserve sub-bundle $C(M)={\mathcal N}^{\flat}$ (see Lemma 1). That
is,
\[
P^{\tau \sigma}\nabla_{\alpha^{\sharp}}\beta_{\sigma} = 0
\]
for any $\beta $ in $C(M)$. Writing this equality in the form
$\left(\nabla_{\alpha^{\sharp}}P\right)^{\tau \sigma}\beta_{\sigma} = 0 $ and using the previous
Corollary we get the result.
\end{proof}
\begin{remark}
Lichnerowicz's  condition, i.e. the requirement that $\nabla P=0$,
guarantees much more than the integrability of the distribution
$\mathcal N$ and, therefore, the local splitting of $M$ into a product
of a symplectic leaf $S$ and complemental manifold $N$ with zero
Poisson tensor. It also guarantees regularity of the Poisson structure,
and reduction of the metric $g$ to the block diagonal form
$g=g_{S}+g_{\mathcal N}$, with the corresponding metrics on the
symplectic leaves and maximal integral manifolds $N_{m}$ of
$\mathcal N$ being independent on the complemental variables (i.e. the metric $g_{S}$ on the symplectic leaves is independent from the
coordinates $y$ along $N_{m}$). Furthermore, the condition $\nabla P=0$ also ensures the independence
of the symplectic form $\omega_{S}$ from the transversal
coordinates $y$ (see Vaisman 1994, Remark 3.11). Finally from
$\nabla ^{g_{S}}\omega_{S}=0$ follows the existence of a
$g_{S}$-parallel Kahler metric on the symplectic leaves.
\end{remark}
\begin{corollary}
Let $\nabla_{\lambda}\omega^{i}=0$ for all $\lambda ,i$ (i.e.
the 1-forms $\omega^{i}=dc^{i}$ are $\nabla^{g}$-covariantly
constant). Then
\begin{enumerate}
\item[i)] The distribution $\mathcal N$ is integrable,
\item[ii)] the vector fields $\xi_{i}$ are Killing vector fields of the metric
$g$, and
\item[iii)] the Casimir functions $c^{i}$ are harmonic: $\Delta_{g}c^{i}=0$.
\end{enumerate}
\end{corollary}
\begin{proof}
The first statement is a special case of (3) in the Theorem above. To prove the second, we
calculate the Lie derivative of $g$ in terms of the covariant derivative $\nabla \omega^{i}$,
\begin{eqnarray}
\nonumber (\mathfrak{L}_{\xi_{i}}g)_{\sigma \lambda}& = & g_{\gamma \lambda}\nabla_{\sigma}\xi^{\gamma}_{i}+ g_{\sigma \gamma}\nabla_{\lambda}\xi^{\gamma}_{i}\\
\nonumber & = & \nabla_{\sigma}\omega^{i}_{\lambda} +  \nabla_{\lambda}\omega^{i}_{\sigma}\\
\nonumber & = & \frac{\partial \omega^{i}_{\lambda}}{\partial x^{\sigma}} + \frac{\partial \omega^{i}_{\sigma}}{\partial x^{\lambda}} - \omega^{i}_{\gamma} (\Gamma^{\gamma}_{\sigma \lambda} + \Gamma^{\gamma}_{\sigma \lambda})\\
\nonumber & = & \frac{\partial^{2}c^{i}}{\partial x^{\sigma}x^{\lambda}} + \frac{\partial c^{i}}{\partial x^{\lambda}x^{\sigma}} - 2\omega^{i}_{\gamma} \Gamma^{\gamma}_{\sigma \lambda}\\
\nonumber & = & 2\frac{\partial^{2}c^{i}}{\partial x^{\sigma}x^{\lambda}} -2\omega^{i}_{\gamma}\Gamma^{\gamma}_{\sigma \lambda}\\
\nonumber & = & 2\nabla_{\lambda}\omega^{i}_{\sigma}.
\end{eqnarray}
Thus, if the condition of the Corollary is fulfilled, $\xi_{i}$ are Killing vector fields. The third statement follows from
\[
\Delta_{g}c^{i}=div_{g}(\xi_{i}=(dc^{i})^{\sharp})=\frac{1}{2}Tr_{g}({\mathcal
L}_{\xi_{i}}g)=\frac{1}{2}g^{\lambda \mu }({\mathcal
L}_{\xi_{i}}g)_{\lambda \mu }.
\]
\end{proof}
\vskip0.5cm

\subsection{Nijenhuis Tensor}
Conventionally the integrability of different geometrical structures presented by a $(1,1)$-tensor
field can be characterized in terms of the corresponding Nijenhuis tensor. Thus, it is interesting
to see the relation of our criteria presented above to the nullity of the corresponding Nijenhuis
tensor.
\begin{definition} Given any $(1,1)$ tensor field $J$ on $M$, there exists a tensor field
 $N_{J}$ of type $(1,2)$ (called the Nijenhuis torsion of $J$) defined as follows
 (see DeLeon and Rodrigues 1989, Sec.1.10):
                \[N_{J}(\xi ,\eta) = [J\xi,J\eta] - J[J\xi,\eta] - J[\xi,J\eta] + J^{2}[\xi,\eta]
                \]
for all vector fields $\xi,\eta $.
\end{definition}

If $J$ is an \textbf{almost product structure}, i.e. $J^{2} = Id$, then $N_{J}=0$ is equivalent
to the integrability of $J$. In fact, given such a structure on $M$, we can define projectors
$v=(1/2)(Id + J)$ and $h=(1/2)(Id -J)$ onto complementary distributions $Im(v)$ and $Im(h)$ in $TM$
such that at each point $x \in M$,
                       \[T_{x}M = Im(v)_{x} \oplus Im(h)_{x}.\]
It is known (see DeLeon and Rodrigues 1989, Sec.3.1) that $J$ is
integrable if and only if $Im(v)$ and $Im(h)$ are integrable, and
that the following equivalences hold:
           \[N_{J} = 0 \leftrightarrow N_{h} = 0 \leftrightarrow N_{v}=0.\]

Consider now the two complementary distributions $T(\mathcal S)$ and $\mathcal{N}$ discussed above. Suppose
that $v$ is $g$-orthogonal projection onto the distribution $T(\mathcal S)$, and $h$ is $g$-orthogonal
projection onto $\mathcal{N}$. Applying these results in this setting we see that that the
distribution $\mathcal{N}$ is integrable if and only if $N_{v}=0$.

Since $v^{2}=v$, and since any $\xi \in T(M)$ can be expressed as $\xi=v\xi + h\xi$, we have
\begin{eqnarray*}
N_{v}(\xi,\eta)
            &=&  [v\xi,v\eta] - v[v\xi,v\eta + h\eta] - v[v\xi+h\gamma,v\eta] + v[v\gamma+h\gamma,v\eta+h\eta]  \\
            &=&  (Id - v)[v\xi,v\eta] + v[h\gamma,h\eta] \\
            &=&  h[v\xi,v\eta] + v[h\xi,h\eta]
\end{eqnarray*}
for all $\xi$ and $\eta$ in $T(M)$. Since $T( \mathcal S)$ is
integrable we have $[v\xi,v\eta]\in T( \mathcal S)$, and so
\[
N_{v}(\xi,\eta) = v[h\xi,h\eta].
\]
As a result, we can restrict $\xi$ and $\eta$ to be sections of
the distribution $\mathcal{N}$ to get the following integrability
condition for $\mathcal{N}$ in terms of (1,1)-tensor $v$:
\[ \mathcal{N}\;  \textrm{is integrable} \; \leftrightarrow \; N_{v}(\xi,\eta)^{\mu}  =
  v^{\mu}_{\nu}[\xi,\eta]^{\nu} \stackrel{*}{=} -\partial_{j}v^{\mu}_{\nu}(\xi \wedge \eta)^{j\nu}=0,\]
$\text{for all}\  \mu\  \text{and all}\ \xi, \eta \in \Gamma(\mathcal{N})$. The equality ($*$) on the right is proved in the
same way as the similar result for the action of $P^{\tau \sigma}$ in the
proof of statement (3) of Theorem 1.
\par
 The tensor $A^{\mu }_{\nu} = g_{\nu \sigma}P^{\sigma \mu }$ discussed
above can be considered to be a linear mapping from $T(M)$ to
$T(\mathcal S)$, but since $A$ is not idempotent, it does not define a
projection. However, the tensors $A$ and $v$, having the same
kernel and image are related in the sense that the
integrability of $\mathcal{N}$ is also equivalent to
\[ A^{\mu}_{\nu}[\xi,\eta]^{\nu} = -\partial_{\sigma}A^{\mu}_{\nu}(\xi \wedge \eta)^{\sigma \nu}=0,\]
for all sections $\xi$ and $\eta$ of the distribution $\mathcal{N}$
(using the same argument as for the tensor $v^{\mu }_{\sigma}$
above).\par

 In fact, since the linear mappings $A,v$ of $T_{m}(M)$ have the same kernel and image for all $m\in M$,  there exists a (non-unique) pure gauge automorphism $D:T(M) \rightarrow T(M)$ of the tangent
bundle (i.e. inducing the identity mapping of the base $M$ and,
therefore, defined by a smooth (1,1)-tensor field $D^{\mu}_{\nu}$)
such that $A^{\mu}_{\sigma } = D^{\mu}_{\nu}v^{\nu}_{\sigma}$. For any
couple $\xi$ and $\eta$ of sections from  $\Gamma ({\mathcal N})$, we have
\begin{eqnarray*}
 A^{\mu}_{\sigma }[\xi,\eta]^{\sigma} & = & \hspace{-2.5mm} -\partial_{\nu}A^{\mu}_{\sigma }(\xi \wedge \eta)^{\nu \sigma } \\
     & = & \hspace{-2.5mm} -\partial_{\nu }D^{\mu }_{\kappa}(v^{\kappa}_{\sigma }(\xi \wedge \eta)^{\nu \sigma} +
     D^{\mu}_{\kappa}\partial_{\nu }v^{\kappa }_{\sigma }(\xi \wedge \eta)^{\nu \sigma}\\
     & = & \hspace{-2.5mm} -D^{\mu }_{\kappa}\partial_{\nu }v^{\kappa}_{\sigma }(\xi \wedge \eta)^{\nu \sigma}\\
                    & = & D^{\mu }_{\kappa}N_{v}(\xi,\eta).
\end{eqnarray*}
This proves
\begin{theorem} There exists a (not unique) invertible linear automorphism $D$ of the bundle $T(M)$
such that for all couples of vector fields $\xi ,\eta \in
\Gamma({\mathcal N})$
\[
A[\xi ,\eta ]=D( N_{v}(\xi ,\eta )).
\]
Thus, $N_{v}\vert_{{\mathcal N}\times {\mathcal N}}\equiv 0$ iff
$A[\xi , \eta]=0\  {\text for\ all}\  \xi , \eta \in \Gamma
({\mathcal N}).$
\end{theorem}

\vskip 1cm

\section{Local criteria for integrability}
Since $M$ is regular, any point in $M$ has a neighborhood in which the Poisson
tensor $P$ has, in Darboux-Weinstein (DW) coordinates
$(y^{A},x^{i})$, the following canonical form (see Weinstein 1983)
\[
P=\begin{pmatrix} 0_{p\times p} & 0_{p \times 2k}\\
0_{2k \times p} & \begin{pmatrix} 0_{k} & -I_{k}\\
I_{k} & 0_{k}
\end{pmatrix}
\end{pmatrix}
\]

We will use Greek indices $\lambda , \mu, \tau $ for general local
coordinates, small Latin $i,j,k$ for the canonical coordinates
along symplectic leaves and capital Latin indices $A,B,C$ for
transversal coordinates. In these DW-coordinates we have, since
$P$ is constant,
\[ (\nabla_{\lambda }P)^{\tau \sigma }=P^{j\sigma }\Gamma^{\tau }_{j\lambda }- P^{j\tau
}\Gamma^{\sigma}_{j\lambda}.
 \]
Using the structure of the Poisson tensor we get, in matrix form,
\[
 (\nabla_{\lambda }P)^{\tau \sigma }=\begin{pmatrix} 0_{p\times p} & P^{js}\Gamma^{T}_{j\lambda}\\
-P^{it}\Gamma^{s}_{j\lambda} & P^{js }\Gamma^{t}_{j\lambda }- P^{jt}\Gamma^{s}_{j\lambda}
\end{pmatrix} ,
\]
where the index $\tau $ takes values $(T,t)$, and the index $\sigma $ takes values $(S,s)$, transversally
 and along the symplectic
leaf respectively.
\par In DW-coordinates we choose $\omega^{\tau}=dy^{\tau}$ as a basis for the co-distribution $C(M)$.
 Now we calculate
(using the symmetry of the Levi-Civita connection $\Gamma$)
\[
 (\nabla_{\lambda }P)^{\tau \sigma }(dy^{I}\wedge dy^{J})_{\alpha \sigma}=
 -\delta^{I}_{\alpha}P^{j\tau}\Gamma^{J}_{j\lambda} +\delta^{J}_{\alpha}
 P^{j\tau}\Gamma^{I}_{j\lambda} ,
\]
so that
\[
 g^{\lambda \alpha}(\nabla_{\lambda }P)^{\tau \sigma }(dy^{I}\wedge dy^{J})_{\alpha \sigma}=
 P^{j\tau}[g^{J\lambda}\Gamma^{J}_{j\lambda} -
 g^{I\lambda}\Gamma^{I}_{j\lambda}].
\]
This expression is zero if $\tau =T$, so the summation goes by $\tau =t$ only.\par

Substituting the Poisson Tensor in its canonical form we get the
integrability criteria (3) of Theorem 1 in the form
\[
g^{J\lambda}\Gamma^{I}_{\lambda t} -
 g^{I\lambda}\Gamma^{J}_{\lambda t}=0,\  \forall \ I,J,t.
\]
Using the metric $g$ to lower indices, we finish the proof of the following
\begin{theorem}
Let $(y^{I},x^{i})$ be local DW-coordinates in $M$. Use capital Latin indices for transversal
coordinates $y$ along $\mathcal N$ and small Latin indices for coordinates $x$ along symplectic leaves. Then
 the distribution $\mathcal N$ is integrable if and only if
\begin{equation}
\Gamma_{JIt}=\Gamma_{IJt},\  \forall \ I,J,t.
\end{equation}
\end{theorem}

\section{Examples: Non-integrability}

\subsection{A Model 4d system}
We now consider a (local) model example of the lowest possible
dimension where the distribution ${\mathcal N}_{g}$
 may not be integrable. This is the case of a 4-d Poisson manifold $(M=\mathbf{R}^{4},P)$ where $rank(P)=2$ at all
  points of the manifold $M$.\par
  Let $P^{ij}$ be the canonical Poisson tensor given in the
 {\bf global} coordinates $x^{\alpha}$ by the following
   $4\times 4$ matrix:
 \[ P=  \begin{pmatrix}
            0 & 0 &0 & 0\\ 0&0&0&0\\0&0& 0 & 1 \\0&0& -1 & 0 \\  \end{pmatrix}
                       \]
Let $\omega^{1}=dx^{1}$, and $\omega^{2}= dx^{2}$. Then $\{\omega^{1}, \omega^{2} \}$ is a basis for
 the kernel $C(M)$ of $P$, and
\begin{equation}\label{C} (\omega^{1} \wedge \omega^{2})_{\alpha \sigma} =
\left\{ \begin{array}{rl}  1, & \alpha = 1 , \, \, \sigma = 2 \\
                          -1, & \alpha = 2 , \, \, \sigma = 1 \\
                           0, & \textrm{otherwise.}   \end{array}\right. \end{equation}

\par
If $h$ is the Euclidian metric in $\mathbf{R}^4$, then it is obvious that the
$h$-orthogonal distribution $\mathcal{N}_{h}$ is generated by the
basic vector fields $\frac{\partial}{\partial x^{1}},
\frac{\partial}{\partial x^{2}}$ and is trivially integrable.\par

Let now $g$ be an arbitrary pseudo-Riemannian metric defined on
$M=\mathbf{R}^4$ by a non-degenerate symmetric (0,2)-tensor $g_{\lambda
\mu}$. The corresponding $g$-orthogonal distribution is denoted by
$\mathcal N$ and the Levi-Civita connection of the metric $g$ by
$\nabla$.\par
 Consider $\nabla_{\lambda}P^{\tau \sigma} =
\partial_{\lambda}P^{\tau \sigma} +
 P^{\tau \mu}\Gamma_{\lambda \mu}^{\sigma} + P^{\sigma \mu}\Gamma_{\lambda \mu}^{\tau}$.
 Since $P$ is constant, the first term of this expression is always zero. Furthermore,
 since each $\omega^{k}$ is in the kernel of $P$, we see that the third term in this expression
 will contract to zero with \ $(\omega^{1} \wedge \omega^{2})_{\alpha \sigma}$. Therefore,
\begin{eqnarray} \nonumber g^{\lambda \alpha}\nabla_{\lambda}P^{\tau \sigma}(\omega^{1}
 \wedge \omega^{2})_{\alpha \sigma} & = & g^{\lambda \alpha}P^{\tau \mu}
 \Gamma^{\sigma}_{\lambda \mu}(\omega^{1} \wedge \omega^{2})_{\alpha \sigma}.\\
\nonumber  & = & g^{\lambda 1}P^{\tau \mu}\Gamma^{2}_{\lambda \mu} - g^{\lambda 2} P^{\tau
\mu}\Gamma^{1}_{\lambda \mu},\hspace{8mm} \mbox{by (\ref{C})}.
\end{eqnarray}
The only values of $\tau $ for which $P^{\tau \mu} \neq 0$ are $\tau = 3$ and $\tau = 4$.
We consider each case individually:\\

\noindent {\boldmath{$\tau = 3$:}} \begin{eqnarray}
\nonumber g^{\lambda \alpha}\nabla_{\lambda}P^{\tau \sigma}(\omega^{1}
 \wedge \omega^{2})_{\alpha \sigma} & = & g^{\lambda 1}P^{3 4}\Gamma^{2}_{\lambda 4} -
 g^{\lambda 2}P^{3 4}\Gamma^{1}_{\lambda 4},\\
\nonumber & = & g^{\lambda 1}\Gamma^{2}_{\lambda 4} - g^{\lambda 2}\Gamma^{1}_{\lambda 4},\\
\nonumber & = & \frac{1}{2}(g^{\lambda 1}g^{2 \delta} - g^{\lambda 2}g^{1 \delta})(g_{\lambda \delta
,4}
 + g_{4 \delta ,\lambda} - g_{\lambda 4,\delta}),\\
\nonumber & = &g^{\lambda 1}g^{2 \delta}(g_{4\delta ,\lambda} - g_{4\lambda ,\delta}).
\end{eqnarray} \\
\noindent {\boldmath{$\tau = 4$:}}\begin{eqnarray}
\nonumber g^{\lambda \alpha}\nabla_{\lambda}P^{\tau \sigma}(\omega^{1} \wedge \omega^{2})_{\alpha
\sigma} &
= & g^{\lambda 1}P^{4 3}\Gamma^{2}_{\lambda 3} - g^{\lambda 2}P^{4 3}\Gamma^{1}_{\lambda 3},\\
\nonumber & = & -g^{\lambda 1}\Gamma^{2}_{\lambda 3} + g^{\lambda 2}\Gamma^{1}_{\lambda 3},\\
\nonumber & = & \frac{1}{2}(-g^{\lambda 1}g^{2 \delta} + g^{\lambda 2}g^{1 \delta})(g_{\lambda
\delta ,3}
 + g_{3 \delta ,\lambda} - g_{\lambda 3,\delta}),\\
\nonumber & = &g^{\lambda 1}g^{2 \delta}(g_{3\lambda ,\delta} - g_{3\delta ,\lambda}).
\end{eqnarray}
\\
Thus, the integrability condition takes the form of the following system of equations
\begin{eqnarray*}
g^{\lambda 1}g^{2 \delta}(g_{3\lambda ,\delta} - g_{3\delta
,\lambda})& =& 0,\\
g^{\lambda 1}g^{2 \delta}(g_{4\delta ,\lambda} - g_{4\lambda ,\delta}) & =& 0
\end{eqnarray*}
equivalent to the symmetry conditions (7).\par
 Clearly both expressions are zero if $g$ is diagonal.
In fact, if $g$ is block-diagonal, then both of the above terms
will also vanish. For these special types of metric, the
transversal
 distribution $\mathcal{N}$ is integrable. For more general metrics, however,  $\mathcal{N}$ may not
 be integrable. For example, let
\[g = \begin{pmatrix} 1 & 0 & f & 0 \\
                      0 & 1 & 0 & 0 \\
                      f & 0 & 1 & 0 \\
                      0 & 0 & 0 & 1 \end{pmatrix},\]
where $f(x)$ satisfies to the condition $\partial_{2}f \neq 0 $.
This symmetrical matrix has $1,1,1+f,1-f$ as its eigenvalues. Thus
$g$ determines the Riemannian metric in the region $|f| < 1$, \,
and the condition
\[g^{\lambda \alpha}\nabla_{\lambda}P^{\tau \sigma}(\omega^{1} \wedge \omega^{2})_{\alpha \sigma} = 0 \]
fails since, for $\tau = 4$ we have:
\begin{eqnarray}
\nonumber g^{\lambda \alpha}\nabla_{\lambda}P^{4 \sigma}(\omega^{1} \wedge \omega^{2})_{\alpha
\sigma} & =
 & g^{\lambda 1}g^{2 \delta}(g_{3\lambda ,\delta} - g_{3\delta ,\lambda}),\\
\nonumber & = & g^{\lambda 1}(g_{3 \lambda ,2} - g_{32,\lambda}),\\
\nonumber & = & g^{\lambda 1}g_{3 \lambda , 2},\\
\nonumber & = & g^{11}g_{31,2} ,\\
\nonumber & = & \partial_{2}f \neq 0.
\end{eqnarray}
\par
As an example of such a function $f$ for which both conditions
(i.e. conditions $\vert f\vert <1$ and $\partial_{2}f \neq 0 $)
are fulfilled in the whole space $\mathbf{R}^4$  we can take the
function $f(x^{1},\ldots ,x^{4})=\frac{1}{\pi}tan^{-1}(x^{2})$
where the principal branch of $tan^{-1}(x)$ is chosen (taking
values between $-\pi/2$ and $\pi/2 $).\par

We can also see that the
distribution $\mathcal{N}$ is not integrable by a direct
computation. Observe that the local basis vectors for
$\mathcal{N}$ are:
\[ \xi_{1} = \partial_{1} +
f\partial_{3}, \,\,\,\, \xi_{2} = \partial_{2}.\]
Their Lie bracket is $[\xi_{1},\xi_{2}] =
\partial_{2}f\partial_{3}$, which is not in the span of
 $\{\xi_{1},\xi_{2}\}$ (since $\partial_{2}f \neq 0$), hence $\mathcal{N}$ is not integrable.
 \vskip 0.5cm

\subsection{Case of a symplectic fibration}
Here we discuss a situation that demonstrates that the very
possibility to choose a global metric $g$ such that the
distribution ${\mathcal N}_{g}$ is integrable is determined mostly
by the topological properties of the
 ``bundle" of leaves of the symplectic foliation, i.e. the existence of a zero curvature Ehresmann connection.\par
 A topologically simple (with regard to the transversal structure) example of a regular Poisson manifold
is a symplectic fibration: a fiber bundle $(M,\pi ,B)$ such that
every fiber $S_{b}=\pi^{-1}(b)$ is endowed with a symplectic
structure $\omega_{S}=\omega_{b}$. Each fiber is symplectomorphic
to the model symplectic manifold $(F,\omega)$, and the transition
functions of a trivialization of this bundle are symplectic
isomorphisms on the fibers (see Guillemin et al 1996, Ch.1). The
inverse $P_{b}=\omega_{b}^{-1}$ of the symplectic form on each
symplectic fiber defines, via the embedding $ \bigwedge
T(S_{b})\rightarrow \bigwedge T(M)$, a smooth (2,0)-tensor field
$P$, i.e. a regular Poisson structure on $M$.
\par
The distribution ${\mathcal N}_{g}$, which is $g$-orthogonal to the fibers
$S_{b}$ with respect to some (pseudo-)Riemannian metric on $M$
(with the condition that the restriction of $g$ to the fiber $S_{b}$
is non-degenerate), defines an Ehresmann connection $\Gamma_{g}$
on the bundle $(M,\pi ,B)$. Integrability of the distribution
${\mathcal N}_{g}$ (i.e. integrability of the connection
$\Gamma_{g}$) means that the curvature (Frobenius tensor) of the
connection $\Gamma_{g}$ is zero. \par

On the bundle $(M,\pi ,B)$ of a symplectic fibration there is a
special class of {\it symplectic connections} $\Gamma $
distinguished by the condition that the holonomy mappings are symplectic diffeomorphisms of the fibers.  It is
proved in Guillemin et al 1996, Ch. I that if $F$ is compact,
connected and simply connected, then for such a connection there
exists a closed 2-form $\omega_{\Gamma}$ on $M$ whose restrictions
to any fiber $S_{b}$ coincide with $\omega_{b}$, and such that the
orthogonal complement $\omega_{\Gamma}$ of the tangent space
$T_{m}(\mathcal S)$ to the fiber passing through a point $m\in M$
is exactly the horizontal subspace $Hor_{\Gamma}(m)$ of the
connection $\Gamma $ at the point $m$. The curvature of the
connection $\Gamma$, which measures the degree of
``non-integrability" of the distribution $Hor_{\Gamma}$, is
determined by the form $\omega_{\Gamma}$ through the {\it
curvature identity} proved in Guillemin et al 1996, Ch. I. Namely,
let $v_{1},v_{2}$ be two arbitrary vector fields on $B$ and denote
by $v^{\sharp}_{1},v^{\sharp}_{2}$ their horizontal lifts to vector fields in $M$.  Then the curvature of
$\Gamma $ is the  vertical (i.e. restricted to the fiber) part
of the 1-form
$i_{[v^{\sharp}_{1},v^{\sharp}_{2}]}\omega_{\Gamma}$, and one has
the equality
\[
-di_{v^{\sharp}_{1}}i_{v^{\sharp}_{2}}\omega_{\Gamma}
=i_{[v^{\sharp}_{1},v^{\sharp}_{2}]}\omega_{\Gamma} \ mod\ B
\]
where $mod\ B$ means ``restricted to the fiber''. This
restriction is zero if and only if the function
$H=i_{v^{\sharp}_{1}}i_{v^{\sharp}_{2}}\omega_{\Gamma}$ is {\bf
constant along the fibers} $S_{b}$.  Then $H=\pi^{*}h$ for some
$h\in C^{\infty}(B)$ is a Casimir function for the Poisson
structure on $M$ constructed as described above.
\par
Having a connection $\Gamma $ (symplectic or not) available on the
bundle $(M,\pi ,B)$, one can define a whole class of
(pseudo-)Riemannian metrics for which the orthogonal complement of
$T(\mathcal S)$ will coincide with $Hor_{\Gamma}$. Namely, we take
a metric $g_{S,b}$ on $T_{m}({\mathcal S})$ smoothly depending on
the point $b$. Then we take an arbitrary metric $g_{B}$ on the
base $B$ and lift it to the horizontal subspaces of $\Gamma$. The
metric $g$ on the total space of the bundle $M$ is now defined by
the condition of orthogonality of $T(\mathcal S)$ and
$Hor_{\Gamma}$. The projection $\pi :M\rightarrow B$ becomes a
(pseudo-)Riemannian submersion (see Cheeger and Ebin 1975).  There
is a relation between the curvatures of $g_{B},g$ and the
curvature of the connection $\Gamma$ (O'Neill formula, see Cheeger and
Ebin 1975, 3.20). Let $m\in M,\ b=\pi(m)$, $X,Y\in T_{b}(B)$ be
two arbitrary tangent vectors at $b$, and let ${\bar X},{\bar Y}\in
Hor_{\Gamma} (m)\subset T_{m}(M)$ be their horizontal lifts at the
point $m$.  Then for the sectional curvatures $K$ of the metric
$g_{B}$ and $\bar K$ of $g$ one has
\[
K_{b}(X,Y)={\bar K}_{m}({\bar X},{\bar Y})+\frac{3}{4}\Vert [{\bar
X},{\bar Y}]^{vert}(m)\Vert ^{2}_{g_{S}\ m}\ ,
\]
where $^{vert}$ means taking the vertical component of the bracket
of horizontal lifts to a neighborhood of $m$ of arbitrary vector
fields in $B$ having values $X,Y$ at the point $b$. Thus, the
curvature of the connection $\Gamma $ measures the difference
between the sectional curvatures of the metric $g_{B}$ and its
horizontal lift to the distribution $Hor_{\Gamma}$.
\par
It is easy to construct examples of bundles which do not allow
integrable Ehresmann connections using the following arguments.
Suppose that a bundle $(M,\pi ,B)$ with a simply-connected base
$B$ allows an integrable Ehresmann connection $\Gamma $. The
holonomy group of the connection $\Gamma $ is discrete (by the
Ambrose-Singer Theorem, since the curvature is zero) and,
therefore, any maximal integral submanifold (say, $V$) of $\Gamma
$ is a covering of $B$. Since $B$ is simply-connected, the
projection $\pi :V\rightarrow B$ is a diffeomorphism. Pick a point
$b\in B$. Then every maximal integral manifold intersects the
fiber $F_{b}$ at one point, defining in this way a smooth
diffeomorphism $q: M\rightarrow F_{b}\simeq F$ smoothly depending
on $b$. Together with the projection $\pi $ this mapping defines a
trivialization $(\pi ,q): M\rightarrow B\times F$ of the bundle
$(M,\pi ,B)$.
\par
This proves the following
\begin{proposition}
Let $(M,\pi ,B;(F,\omega ))$ be a symplectic fibration with the
model symplectic fiber $(F,\omega)$, base $B$ and the total
Poisson space $M$. If the bundle $(M,\pi ,B)$ is topologically
non-trivial, then the Poisson manifold $(M,P=\omega_{b}^{-1}$
cannot be endowed with a (global) pseudo-Riemannian metric $g$
such that the orthogonal distribution $\mathcal{N}_{g}$ would be
integrable.
\end{proposition}
\par
Thus, if we take an arbitrary nontrivial bundle over a
simply-connected manifold $B$,  it can not have a nonlinear
connection of zero curvature. An example is the tangent bundle
$(T(\mathbf{C}P(2)),\pi ,\mathbf{C}P(2))$ over $B=\mathbf{C}P(2)$,
where the standard symplectic structure on $B=\mathbf{C}P(2)$
determines a (constant) symplectic structure along the fibers.\par

\section{Examples of Integrability: Linear Poisson structure}
Let $\mg$ be a real n-dimensional Lie algebra with a basis
$\{e_{\mu }\}$ and Lie bracket $[e_{\mu },e_{\nu }]=c_{\mu \nu
}^{\sigma }e_{\sigma }$. Let $G$ be a connected Lie group with Lie algebra $\mg$. The Killing form $K$ on $\mg$ is the invariant,
symmetric, bilinear form defined by
\[
K(x,y)=Tr(ad(x)\circ ad(y)),\;\;  K_{\mu \nu } =
Tr(ad_{e_{\mu }}\cdot ad_{e_{\nu }}),
\]
where\, $ad_{v}(X) = [v,X]$,\, $X \in \mg $ (see Bourbaki 1968,
Ch.3). Let $\{f^{\nu }\}$ be the dual basis on the dual space
$\mgs$, and let $\lambda_{\nu }$ be coordinates for $\mgs$
relative to this basis: $\lambda =\lambda_{\nu }f^{\nu },\
\lambda_{\nu }=<\lambda ,e_{\nu }>$. The dual space $\mgs$ with its linear Lie-Poisson structure
\begin{equation}
P^{\mu \nu }(\lambda)= \{\lambda_{\mu },\lambda_{\nu }\}=c_{\mu
\nu }^{\sigma }\lambda_{\sigma },
\end{equation}
is a model example of a Poisson manifold. The adjoint action $Ad(g)$
of the corresponding Lie group $G$ on $\f{g}$ defines the linear
co-adjoint action $Ad^{*}(g)$ of $G$ on $\f{g}^{*}$.
 The symplectic leaves of the Lie-Poisson structure (9) are co-adjoint orbits
of $G$ (see Kirillov, 2004). \par
Casimir functions are exactly the $Ad^{*}(G)$-invariant functions on $\mgs.$
  In many cases (for instance for real semi-simple Lie groups) one can choose $k=rank(G)$
 polynomial Casimir functions $c_{i}$ that are functionally independent on $M= \mgs_{reg}$,
  and any Casimir function is a function of the polynomials $c_{i}$ (see Kirillov, 2004).\par

 The tangent space to $\f{g}^{*}$ at each point
can be identified with the (vector) space $\f{g}^{*}$ itself and,
correspondingly, the tangent bundle $T(\f{g}^{*})$ splits
$Ad^{*}(G)$-equivariantly: $T(\f{g}^{*})\simeq \f{g}^{*}\times
\f{g}^{*}$. The cotangent bundle $T^{*}(\f{g}^{*})$ takes the form
$T^{*}(\f{g}^{*})\simeq \f{g}^{*}\times \f{g}.$\par
  The action of $G$ induced by $Ad^{*}(G)$ on $T(\f{g}^{*})\simeq \f{g}^{*}\times \f{g}^{*}$ is
$Ad^{*}(g)\times Ad^{*}(g)$, while the dual action on the cotangent
bundle $T^{*}(\f{g}^{*})\simeq \f{g}^{*}\times \f{g}$ takes the
form $Ad^{*}(g)\times Ad(g).$  In particular, the action on the second
factor coincides with the adjoint action of $G$ on $\f{g}$.
 Below we will be using these identifications without further comments.\par

If $\mg$ is a real semi-simple Lie algebra, the Killing form $K$
is non-degenerate and can be used to identify $\mg$ with $\mgs$.
Under this identification, the adjoint action of $G$ corresponds
to the co-adjoint and, correspondingly, adjoint orbits correspond
to co-adjoint ones. Thus, one can translate the linear Poisson
structure to the Lie algebra $\mg $ and use the available information
about the (singular) foliation of $\mg$ by the adjoint orbits (see
Warner 1972, Sec.1.3) to study the $K$-orthogonal distribution
${\mathcal N}_{K}$ defined by the (pseudo-Riemannian) Killing
metric $K$.\par

\subsection{Compact semi-simple Lie algebra}
Consider the case when $\mathfrak g$ is a compact real semi-simple
Lie algebra, i.e. the Lie algebra of a compact semi-simple Lie
group. Then the Killing form $K$ is negative definite, and, therefore,
$-K$ is an invariant Riemannian metric on $\mg$ (see Knapp, 1996,
Ch.4).\par

The canonical isomorphism $T_{x}^{*}(\mgs )\simeq \mg $ defined above
allows us to consider the Killing form $-K_{\mu \nu }$ on $\mg$ as
a covariant metric $g^{\mu \nu }$ on $\mgs$. Basic vectors
$e_{\mu }\in \mg $ are identified with the covectors
$d\lambda_{\mu }$ in $\mgs$. To be consistent with the upper/lower
indices duality we will denote these {\it basic covectors} in
$\mgs$ by $e^{\mu }$.
\par

 Recall that we have the following condition (see (3) in Theorem 1) for the $g$-orthogonal space
  $\mathcal{N}$ to be integrable.
\[
\mathcal{N}\,\, \textrm{is integrable}\,\, \Leftrightarrow \,\,
 g^{\gamma \alpha}\nabla_{\gamma}P^{\tau \sigma}(\omega_{\alpha}\eta_{\sigma} - \omega_{\sigma}\eta_{\alpha}) = 0,
  \,\,\,\,\, \forall \tau,
\]
where $\omega$ and $\eta $ are any two elements in the kernel of $P$. Using this condition we can
prove the following
\begin{theorem} Let $\mathfrak g$ be a compact semi-simple Lie algebra and let
 $M={\mathfrak g}^{*}_{reg}$ with the standard linear Lie-Poisson structure.  Let $g$ be the inverse
 to the metric on $\mathfrak g$ given by the restriction of the negative Killing form on $\mathfrak g$.
 Then the distribution $\mathcal{N}$ is
integrable if and only if $C(M)$ is an abelian subalgebra of \,
$T^{*}(M)$ endowed with the (pointwise on ${\mathfrak g}^{*}$)
bracket induced from the Lie algebra $\mg $ via the identification
$T^{*}({\mg}^{*})\simeq {\mg}^{*}\times {\mg}$.
\end{theorem}
\begin{proof} Let $\omega = \omega_{\alpha} e^{\alpha }$ and $\eta = \eta_{\beta} e^{\beta}$ be two (local) sections of $C(M)$. Then we have
\begin{eqnarray}
\nonumber g^{\gamma \alpha}\nabla_{\gamma}P^{\tau \sigma}(\omega_{\alpha}\eta_{\sigma} - \omega_{\sigma}\eta_{\alpha}) & = & g^{\gamma \alpha}\nabla_{\gamma}P^{\tau \sigma}(\omega_{\alpha}\eta_{\sigma} - \omega_{\sigma}\eta_{\alpha}),\\
\nonumber & = & g^{\gamma \alpha} c_{\tau \sigma}^{\gamma}(\omega_{\alpha}\eta_{\sigma} -
\omega_{\sigma}\eta_{\alpha}).
\end{eqnarray}
Here the expression (9) for the Poisson tensor was used as well as
the flatness of the metric $-K$ on $M={\mathfrak g}^{*}_{reg}$.
Observe that
\begin{eqnarray}
\nonumber g^{\gamma \alpha}c_{\tau \sigma}^{\gamma} & = & g^{\gamma \alpha}[e^{\tau},e^{\sigma}]_{\gamma},\\
\nonumber &=& -<e^{\alpha}, [e^{\sigma},e^{\tau}]>_{g},\\
\nonumber &=& -<[e^{\alpha},e^{\sigma}],e^{\tau}>_{g},\\
\nonumber &=& -g^{\gamma \tau}c_{\alpha \sigma}^{\gamma}.
\end{eqnarray}
Hence,
\begin{eqnarray}
\nonumber g^{\gamma \alpha}\nabla_{\gamma}P^{\tau \sigma}(\omega_{\alpha}\eta_{\sigma} - \omega_{\sigma}\eta_{\alpha}) & = & - g^{\gamma \tau}c_{\alpha \sigma}^{\gamma}(\omega_{\alpha}\eta_{\sigma} - \omega_{\sigma}\eta_{\alpha}),\\
\nonumber &=& -g^{\gamma \tau}(c_{\alpha \sigma}^{\gamma} - c_{\sigma \alpha}^{\gamma})\omega_{\alpha}\eta_{\sigma},\\
\nonumber &=& -2g^{\gamma \tau}c_{\alpha \sigma}^{\gamma}\omega_{\alpha}\eta_{\sigma},\\
\nonumber &=& -2g^{ \tau \gamma}[\omega,\eta ]_{\gamma}.
\end{eqnarray}
\end{proof}

Using the Killing form to identify $\mg $ and $\mgs$, we can consider $M$ as
$M=({\mg}_{reg},P).$ Thus, we have the following
\par
\begin{corollary}
The distribution $\mathcal N$ on the manifold ${\mathfrak g}^{*}_{reg}$ for a compact semi-simple
Lie algebra $\mathfrak g$ is {\bf integrable}.  Furthermore, via the identification of \ $\mgs$ with
$\mg$ as above,
 each connected component (Weyl Chamber) of the Lie algebra $\mathfrak{t}$ of a maximal  torus $T\subset G$
  is a maximal integral surface of the distribution $\mathcal N$ at each point $x$.
\end{corollary}
\begin{proof}
Let $\mathfrak t$ be one of the maximal commutative subalgebras of $\mg $ (the Lie algebra of a
maximal torus $T\subset G$).\par

 Recall that the root decomposition delivers the $K$-orthogonal decomposition of $\mg $:
\[
\mg ={\mathfrak t}\oplus \sum_{\alpha \in \Sigma}\mg ^{\alpha},
\]
where $\Sigma $ is the root system of the couple $({\mathfrak g}_{c}, {\mathfrak t}_{c})$ and
$\mg^{\alpha}=\mg \cap {\mathfrak g}_{c}^{\alpha}$.\par

Any connected component (Weyl Chamber) of $\mathfrak
t$ is a maximal integral surface of the distribution $\mathcal N$
at each regular point $x$ since $ \mathfrak t$ is $K$-orthogonal
to the tangent space of each adjoint orbit (symplectic leaf) in
$M={\mathfrak g}^{*}_{reg}$ (see Knapp, 1996, Ch.4)
  \[
  T_{x}(Ad(G)X)=\left\{ x+ \sum_{\alpha \in \Sigma} {\mathfrak g}^{\alpha} \right\}.
  \]

  Through each point $x\in \mg $ there passes at least one such subspace $\mathfrak t$, and
 a point $x$ is regular if and only if this $\mathfrak t$, containing $x$ is unique. This proves the
 statement.
\end{proof}

\subsection{Non-compact semi-simple Lie algebras}
Let $G$ be a connected real semi-simple Lie Group with Lie Algebra $\mg $. As in the compact case, symplectic orbits of
$\mg $ with this structure are exactly adjoint orbits of $G$, but
 in contrast to the compact case, an adjoint orbit of $X\in {\mg}$ is closed iff $X$ is semi-simple, i.e  $ad(X)$  is semi-simple (see Warner 1972, Prop.1.3.5.5).
 The subset ${\mg}'$ of {\bf regular} elements, i.e. semi-simple elements $X$ with  centralizer $X^{\mg}$ of minimal dimension (see Warner 1972, Sec.1.3.4) endowed with the induced Poisson structure is
 an open and dense subset of $\mg$. Its structure is as follows.  Let $\mathfrak j$ be a Cartan subalgebra of $\mg $.
Let ${\mathfrak j}'={\mathfrak j}\cap {\mg}'$. Put
\[
{\mg}({\mathfrak j})=\bigcup_{x\in G}Ad(x){\mathfrak j}',
\]
where $G=Int({\mg})$ is the adjoint group of $\mg$. Then (see Warner, Prop. 1.3.4.1),
\[
{\mg}'=\bigcup_{l}{\mg}({\mathfrak j}_{l}),
\]
where ${\mathfrak j}_{l}$ for $1\leq l\leq q$ are representatives of (a finite number of) conjugacy
classes of Cartan subalgebras of $\mg$.\par

Now, pick $1\leq l\leq q$ and let ${\mg}={\mathfrak k}\oplus \mu $ be the Cartan decomposition of $\mg$
 such that
\[
{\mathfrak j}_{l}={\mathfrak j}_{l\ {\mathfrak k}}\oplus {\mathfrak j}_{l\ \mu}  ={\mathfrak
j}_{l}\cap {\mathfrak k}\oplus {\mathfrak j}_{l}\cap \mu
\]
is the direct sum decomposition of the Cartan subalgebra
${\mathfrak j}_{l}$ into compact and noncompact parts.  It is
known that the Killing form $K$ is positive definite on $\mu $ and
negative definite on $\mathfrak k$. Using the Cartan decomposition
of ${\mathfrak j}_{l}$ above, we see that the restriction of the
Killing form to the subspace $X+{\mathfrak j}_{l}\subset
T_{X}({\mg})$, and all its conjugates, has constant signature and
is non-degenerate at all points $X\in \bigcup_{x\in
G}Ad(x){\mathfrak j}'$. Therefore, the same is true for its
$K$-orthogonal complement. The restricted root decomposition
\[
{\mg}={\mathfrak j}_{l}\oplus \sum_{\alpha \in \Sigma }{\mg}^{\alpha},
\]
where $\Sigma $ is the system of (restricted) roots of the pair $({\mg},{\mathfrak j}_{l\ \mu})$ can
be used to show that $X+{\mathfrak j}_{l}$ is the $K$-orthogonal complement to $T_{X}(Ad(G)X)$ in
$T_{X}({\mg}).$\par

We call an element $\lambda \in {\mgs}$ $*$-regular if the corresponding element
$X_{\lambda}+i_{K}^{-1}\lambda $ of $\mg$ is regular. Then, arguments similar to those in the
previous subsection can be used to prove the following
\begin{proposition}
Let $\mg $ be a real semi-simple Lie algebra. Consider the dual
space $({\mg}^{*},P)$ with its linear Poisson structure. Endow
${\mg}^{*}$ with the (pseudo)-Riemannian metric $K^{*}$ induced by
the Killing form on $\mg$. Let $M={\mg}^{*}_{reg}$ be the (open
and dense) submanifold of ${\mg}^{*}$ of co-adjoint orbits
(symplectic leaves) of $*$-regular elements.  Then the restriction
of $K^{*}$ to each orbit in $M$ has constant signature and is
non-degenerate. The $K^{*}$-orthogonal
 distribution ${\mathcal N}$ to the symplectic leaves is integrable. Maximal integral
submanifolds of $\mathcal N$ are images under the identification $i_{K}:{\mg}\equiv {\mgs}$ of
(the regular parts of) Cartan subalgebras of $\mg$.
\end{proposition}
\begin{remark}
We had to add the condition of semi-simplicity of an element because of the presence in $\mathfrak g$ of
a principal nilpotent orbit of the same maximal dimension. The restriction of the Killing form to
such orbits is degenerate. For example, consider the Lie algebra $\mg =\f{sl}(2,\mathbf{R})$. In addition to the
closed semi-simple adjoint orbits of elliptic and hyperbolic elements, there is also the adjoint
nilpotent orbit in $\mg $ of the same dimension 2.
\end{remark}

\subsection{Dual to the Euclidian Lie algebra $\f{e}(3)$.}

 Let $\f{g}$ be the Lie algebra $\f{se}(3)$ of the group of proper Euclidian motions in $\mathbf{R}^3$, let $\f{g}^{*}$ be its dual, and let $M=\f{g}^{*}_{reg}$ be the (open, connected and dense) subspace of 4d
 co-adjoint orbits in $\f{g}^{*}.$\par
   The Killing form of the Lie algebra $\f{se}(3)$ is degenerate, so we identify $\f{g}^{*}$ with
 $\f{g}$ via the Euclidean scalar product (see Marsden and Ratiu 1994, Ch.8). Elements of $\f{g}$
can be represented as vectors $(\bb x, \bb p)$ in $\f{so}(3)\oplus \bf{R}^{3}$ with the Lie bracket
defined as
\[[(\bf x,p),(x',p')]= (x\times x',\ x\times p' + \ p\times x').\]
\par
We can consider vectors $\bf x$ in $\f{s0}(3)$ to be skew-symmetric matrices due to the isomorphism
\[ \bb x \rightarrow \bhat{x}= \begin{pmatrix} 0 & -x_{3} & x_{2}\\
                                                x_{3} & 0 & -x_{1}\\
                                                -x_{2} & x_{1} & 0 \end{pmatrix}. \]

 The canonical linear Poisson structure on $M\subset {\f g^{*}} \simeq {\f g}$ has, in this notation,
  the Poisson tensor has the form
\[P(\bb x, \bb p) =  \begin{pmatrix} \bhat x & \bhat p \\ -\bhat p & 0 \end{pmatrix}.\]

Casimir functions of this structure are $c_{1}={\bb x}\cdot {\bb
p}$ and $c_{2}={\bb p}\cdot {\bb p}$. The subspace of regular (4d)
co-adjoint orbits is, in this notation, defined by the condition
$c_{2}({\bb x}, {\bb p}) \ne 0$.\par

For any $\bb y$ we have $\bb x \times \bb y = -\bb y \bhat x = \bhat x \bb y^{\rm T}$. This allows
us to express the adjoint action on $\f{se}(3)$ as
\[ad_{(\bb y, \bb q)}(\bb x,\bb p) = -(\bb x \bhat y , \bb p \bhat y + \bb x \bhat q)=
\begin{pmatrix} \bhat y &  0 \\ \bhat q & \bhat y \end{pmatrix}(\bb x ,\bb p )^{\rm{T}}. \]

\par
We want to construct a symmetric $(2,0)$ tensor $g$ that is invariant under the adjoint action of the Lie algebra that we can possibly use as a Riemannian metric. Suppose then, that we have a non-degenerate scalar product $g_{0}$
defined on $\f g$  by a constant symmetric matrix:
\[<(\bb x, \bb p),(\bb x', \bb p')>_{g_{0}} = (\bb x, \bb p)\begin{pmatrix} A & B \\ B^{\rm T} & C \end{pmatrix}(\bb x', \bb p')^{\rm T} .\]
We extend $g_{0}$ to a covariant metric $g$ on $M$ by setting
$g(\bb x , \bb p) = g_{0}$. We would like to choose this metric to
be invariant under the co-adjoint action.  Thus, $g$ should
satisfy the equation (we use the identification of the tangent and
cotangent bundles to $\f{g}^{*}$ as above)

\[
<ad_{(\bb y. \bb q)}(\bb x, \bb p), ( \bb x', \bb p') >_{g} +
<(\bb x, \bb p),ad_{(\bb y, \bb q )}(\bb x', \bb p')>_{g}=0.
\]

Since this must hold for arbitrary vectors $(\bb x,\bb p)$ and $(\bb x',\bb p')$ in $\f g$, we have
the following condition on $g$.
\[ \begin{pmatrix}\bhat y^{\rm T} & \bhat q^{\rm T} \\ 0 & \bhat y ^{\rm T} \end{pmatrix}\cdot g \, +\,  g \cdot  \begin{pmatrix} \bhat y & 0 \\ \bhat q & \bhat y \end{pmatrix} = 0.\]
This condition is equivalent to the following system of equations:
\begin{eqnarray*}
1)& & A\bhat y - \bhat y A + B \bhat q - \bhat q B^{\rm T} = 0\\
2)& & B \bhat y - \bhat y B - \bhat q C = 0\\
3)& & B^{\rm T} \bhat y - \bhat y B^{\rm T} + C \bhat q = 0\\
4)& & C\bhat y - \bhat y C = 0.
\end{eqnarray*}
Since these equations must be valid for arbitrary $\bhat y$ and $\bhat q$, it is easy to see that
the metric $g$ must be of the form

\begin{equation}
 g =  \begin{pmatrix} \alpha I & \beta I \\ \beta
I & 0
\end{pmatrix}.
\end{equation}

Thus, the only $ad$-invariant (constant) metrics on $M$ are those
having this special form. Note that such a metric cannot be
Riemannian (see Zefran et al 1996). In fact, the distinct
eigenvalues of such a metric are $\lambda_{i} = (\alpha \pm
\sqrt{\alpha^{2} + 4\beta^{2}})/2$,\;  $i=1,2$ (of multiplicity 3
each). The product of these eigenvalues is
$\lambda_{1}\lambda_{2}= -\beta^{2}\leq 0$. Thus, if
non-degenerate (i.e. $\beta \neq 0$), the metric $g$ has signature
$(3,3)$.\par

Let $T_{(\bb x ,\bb p )}(\mathcal S)$ be the space tangent to the
symplectic leaf passing through the point $(\bb x ,\bb p ) \in M
$. Then we can define its $g^{-1}$-orthogonal complement
$\mathcal{N}_{(\bb x , \bb p)}$, and a $g^{-1}$-orthogonal
distribution $\mathcal{N}$ on $M$. The covectors $\omega^{1}
=dc_{1}= (\bb p, \bb x)$ and $\omega^{2}= dc_{2}= (\bb 0, \bb p)$
form a basis for the subspace $C_{m}(M) = ker(P_{m})\subset
T^{*}_{m}(M) \equiv \f g$, and the tangent vectors
$\xi_{1}=(\omega^{1})^{\sharp}$ and $\xi_{2} = (\omega^{2})^{\sharp}$ form, at each point $m=(\bb x ,\bb p )$, a basis for
$\mathcal{N}_{m}$. We have
\[
\xi_{1}=\begin{pmatrix} \alpha {\bb p}+\beta {\bb x}\\ \beta {\bb
p}\end{pmatrix},\ \xi_{2}=\begin{pmatrix}
\beta {\bb p}\\
0\end{pmatrix}
\]\par

Consider the case $\alpha=0$. It is easy to see that, in the
basis $\xi_{i}$, the  restriction of the metric $g$ to the
subspace $\mathcal{N}_{m}$ of the  distribution $\mathcal{N}$ at a
point $m=(\bb x ,\bb p )$ has the form
\[
g=\begin{pmatrix} 2c_{1}(\bb x ,\bb p) & c_{2}(\bb x ,\bb p)\\
c_{2}(\bb x ,\bb p) & 0
\end{pmatrix} .
\]
Since $det(g)=-({\bb p}\cdot {\bb p})^{2}$, this restriction is
non-degenerate on the subset of regular (4d) co-adjoint orbits.
On the distribution $\mathcal{N}$ the metric $g$ has signature
(1,1). Thus, on the tangent spaces $T(S)$ of symplectic foliation,
$g$ has signature (2,2) at all (regular) points, and results of
Sec.3 are applicable here.

However, the methods from that section are not necessary in this case. Since we have explicit expressions for the
 vectors $\xi_{i}$, and for the tensors $P$ and $g$, it is easy to check the integrability of
 $\mathcal{N}$ directly. The distribution $\mathcal{N}$ will be integrable if and only if the
 Lie bracket $[\xi_{1},\xi_{2}]$ of vector fields $\xi_{i}$ considered as $\mathfrak g$-values functions on $M$
  remains in $\mathcal{N}$.

\begin{proposition} For any choice of (constant) non-degenerate $ad$-invariant metric $g$ on the subspace
 $M=\f{g}^{*}_{reg}$ of 4d co-adjoint orbits of dual space $e(3)^{*}$ of the 3-d Euclidian lie algebra $e(3)$
 , the distribution $\mathcal{N}$ is integrable. For metrics (8) with $\alpha =0$, the maximal integral
 submanifold passing through a point $({\bb x},{\bb p})$ is
 presented, in parametrical form as
 \[
(s,t)\rightarrow e^{s}\begin{pmatrix}  {\bb x}^{T}\\{\bb p}^{T}\end{pmatrix} +e^{t}\begin{pmatrix}
{\bb p}^{T}\\{\bb 0}
\end{pmatrix}.
 \]
\end{proposition}
\begin{proof} We calculate the Lie bracket of the basis for $\mathcal{N}$.
\begin{align*}
[\xi_{1},\xi_{2}] & = \xi_{2}^{ \rm T}(\xi_{2}^{\rm T}) - \xi_{1}^{\rm T}(\xi_{2}^{\rm T})\\
    & =  \begin{pmatrix} 0 & \beta \\ 0 & 0 \end{pmatrix} \begin{pmatrix} \alpha \bb p^{\rm T} + \beta
    \bb x^{\rm T} \\ \beta \bb p^{\rm T} \end{pmatrix} - \begin{pmatrix} \beta & \alpha \\ 0 & \beta
    \end{pmatrix}\begin{pmatrix} \beta \bb p^{\rm T} \\ 0 \end{pmatrix}\\
    &=0.
\end{align*}
The last statement is easily checked by direct calculation.
\end{proof}
\section{Conclusion}
In this work we discuss necessary and sufficient conditions for
the distribution $\mathcal N$ on a regular Poisson manifold
$(M,P)$ defined as orthogonal complement of tangent to symplectic
leaves with respect to some (pseudo-)Riemannian metric $g$ on $M$
to be integrable. We present these conditions in different forms,
including a condition in terms of a symmetry of Christoffel
coefficients of the Levi-Civita connection of the metric $g$ and get
some corollaries, one of which specifies the part of the Lichnerowicz
($\nabla P=0$) condition ensuring integrability of $\mathcal N$
(see Vaisman 1994, 3.11). We present examples of non-integrable
$\mathcal N$ (the model 4d case and the case of a nontrivial
symplectic fibration). We prove integrability of $\mathcal N$ on
the regular part of the dual space $\mgs$ of a real semi-simple
Lie algebra $\mg$ and the same in the case of the 3d Euclidian Lie
algebra $\f{e}(3)$ with a linear Poisson structure.
\par As the case of a symplectic fibration shows, the
integrability of $\mathcal N$ is possible only on a topologically
trivial bundle (trivial transversal topology). Thus, it would be
interesting to study maximal integral submanifolds of $\mathcal N$
in the case of nontrivial symplectic bundles. In particular, it
would be interesting to get conditions on the metric $g$ under
which these maximal integral submanifolds would have maximal
possible dimension. Probably the methods of Cartan-Kahler theory
(see Bryant et. al., 1991) can be employed to investigate these
questions.\par

The authors would like to express their gratitude to Professor
V. Guillemin for his interest in this work.

\section{References}

\par
1. A.Bloch, P.S.Krishnaprasad, J.Marsden, and T.Ratiu, The
Euler-Poincare Equations and Double Bracket Dissipation, Comm.
Math. Phys. (1996), \textbf{175}, 1-42.\par

2. N.Bourbaki, {\it Groupes et Algebres de Lie}, Hermann, Paris,
1968.\par

3. R.L.Bryant, S.S.Chern, R.B.Gardner, H.L.Goldschmidt,
P.A.Griffiths, {\it Exterior Differential Systems}, MSRI Publ.,
v.18, Springer-Verlag, N.Y., 1991.

4. J.Cheeger, D.Ebin, {\it Comparison Theorems in Riemannian
geometry}, North-Holland, Amsterdam, 1975.\par

5. M. de Le\'{o}n, and P.Rodrigues, {\it Methods of Differential
Geometry in Analytical Mechanics}, Elsevier, Amsterdam, 1989.\par

6. M.Grmela, Common Structure in the Hierarchy of Models of
Physical Systems, in "Continuum Models and Discrete Systems", v.1,
1990.\par

7. V.Guillemin, , E.Lerman, and S. Sternberg, {\it Symplectic
Fibrations and multiplicity diagrams}, Cambridge Univ. Press,
Cambridge, 1996.\par

8. D.Fish, {\it Geometrical structures of the metriplectic
dynamical systems}, Ph.D Thesis, Portland, 2005.\par

9. A.Kirillov, {\it Lectures on the Orbit Method},
Springer-Verlag, N.Y., 2004.

 10. A.Knapp,  {\it Lie Groups Beyond an
Introduction}, Birkh\"{a}user, Boston, 1996.\par

11. P.Libermann, Pfaffian systems and Transverse Differential
Geometry, in Cahen, Flato (ed) "Differential Geometry and
Relativity", pp. 107-126, D.Reidel, Dordrecht, 1976.\par

12.  J.Marsden, and T.Ratiu, {\it Introduction to Mechanics and
Symmetry}, Springer-Verlag, NY, 1994.\par

13. P.Molino, { \it Riemannian foliations}, Birkh\"{a}user Verlag,
Basel, 1988.\par

14. P.Morrison, {A Paradigm for Joined Hamiltonian and Dissipative
Systems}, Physica \textbf{18D}(1986), 410-419.\par

15.  B.O'Neill,  {\it Semi-Riemannian Geometry}, Academic Press,
1983.\par

16. B.Reinhart,  {\it Differential Geometry of Foliations},
Springer, 1983.\par

17. I.Vaisman, {\it Lectures on the Geometry of Poisson
Manifolds}, Birkh\"{a}user Verlag, Basel, 1994.\par

18.  G. Warner,  { \it Harmonic Analysis on Semi-simple Lie
groups}, vol.1, Springer,1972.\par

19. A. Weinstein, {The local structure of Poisson manifolds},
J.Differential Geometry, v.18, 1993, pp.523-557.\par

20.  M.Zefran, V.Kumar, and C.Croke, {Choice of Riemannian Metrics
for Rigid Body Kinematics}, Proc. 1996 ASME Design Eng. Tech.
Conf., 1996.

\end{document}